\documentclass[12pt]{amsart} 
\usepackage{marvosym}
\usepackage{amsmath}
\usepackage{amssymb}
\usepackage{amsthm}
\usepackage{color}
\usepackage{graphicx}
\usepackage[applemac]{inputenc}
\usepackage{phonetic}
\usepackage{pb-diagram}
\usepackage[T1]{fontenc}
\usepackage{ stmaryrd }
\usepackage{pst-all}
\usepackage{url}
\newtheorem{prop}{Proposition}[section]
\newtheorem{definicion}[prop]{Definition}
\newtheorem{lema}[prop]{Lemma}
\newtheorem{teor}[prop]{Theorem}
\newtheorem{hecho}[prop]{Fact}
\newtheorem{obs}[prop]{Observation}
\newtheorem{preg}[prop]{Question}
\newtheorem{corol}[prop]{Corollary}
\newtheorem{conj}[prop]{Conjecture}
\newtheorem{const}[prop]{Construction}
\theoremstyle{definition}
\newtheorem{ejem}[prop]{Example}
\newtheorem{rmk}[prop]{Remark}

\def\Ind#1#2{#1\setbox0=\hbox{$#1x$}\kern\wd0\hbox to 0pt{\hss$#1\mid$\hss}     
\lower.9\ht0\hbox to 0pt{\hss$#1\smile$\hss}\kern\wd0}
\def\Notind#1#2{#1\setbox0=\hbox{$#1x$}\kern\wd0\hbox to 0pt{\mathchardef
\nn=12854\hss$#1\nn$\kern1.4\wd0\hss}\hbox to
0pt{\hss$#1\mid$\hss}\lower.9\ht0 \hbox to
0pt{\hss$#1\smile$\hss}\kern\wd0}

\newcommand{\bp}{\begin{prop}}
\newcommand{\ep}{\end{prop}}
\newcommand{\bd}{\begin{definicion}}
\newcommand{\ed}{\end{definicion}}
\newcommand{\bej}{\begin{ejem}}
\newcommand{\eej}{\end{ejem}}
\newcommand{\bl}{\begin{lema}}
\newcommand{\el}{\end{lema}}
\newcommand{\bh}{\begin{hecho}}
\newcommand{\eh}{\end{hecho}}
\newcommand{\bpreg}{\begin{preg}}
\newcommand{\epreg}{\end{preg}}
\newcommand{\bo}{\begin{obs}}
\newcommand{\eo}{\end{obs}}
\newcommand{\bcon}{\begin{conj}}
\newcommand{\econ}{\end{conj}}
\newcommand{\brmk}{\begin{rmk}}
\newcommand{\ermk}{\end{rmk}}
\newcommand{\bc}{\begin{corol}}
\newcommand{\ec}{\end{corol}}
\newcommand{\bconst}{\begin{const}}
\newcommand{\econst}{\end{const}}
\newcommand{\bdem}{\begin{proof}}
\newcommand{\edem}{\end{proof}}
\newcommand{\benum}{\begin{enumerate}}
\newcommand{\eenum}{\end{enumerate}}
\newcommand{\bitem}{\begin{itemize}}
\newcommand{\eitem}{\end{itemize}}
\newcommand{\be}{\begin{ejem}}
\newcommand{\ee}{\end{ejem}}
\newcommand{\bt}{\begin{teor}}
\newcommand{\et}{\end{teor}}

\newcommand{\R}{\mathbb{R}}

\newcommand{\Z}{\mathbb{Z}}

\newcommand{\C}{\mathfrak{C}}

\newcommand{\Nat}{\mathbb{N}}

\newcommand{\ov}[1]{\overline{#1}}
\newcommand{\Le}{\mathcal{L}}
\title{A note on pseudofinite dimension and forking}

 \author{Darío Alejandro García}
\address{Darìo Alejandro Garc\'ia  \\
Departamento de Matem\'aticas \\
Universidad de los Andes \\
Cra 1 No. 18A-10, Edificio H \\
Bogot\'a, 111711 \\
Colombia}

\email{dagarciar@gmail.com}

\begin{document}
\maketitle

\textbf{Abstract.} In this paper we show that an instance of dividing in pseudofinite structures can be witnessed by a drop of the pseudofinite dimension. As an application of this result we give new proofs of known results for asymptotic classes of finite structures.
\section{Introduction}
Dimension theory (so-called) has become one of the most important concepts in model theory and has been used to give a combinatorial description of the definable sets of first order structures.  Even more, it is possible to get structural properties of the models of a first order theory $T$ by assuming some bound on the different ranks associated to $T$. 

One of the recurrent themes in the notions of rank is their relationship with forking. It is often desired that any instance of forking (on types or formulas) can be detected by a decrease of the dimension.

In \cite{HW}, Hrushovski and Wagner defined the notion of quasidimension on some structure $M$ as a way to generalize the concept of dimension allowing values different from the integers. The main example is what I call ``logarithmic pseudofinite dimension'' which is  defined on ultraproducts of finite structures by taking the logarithm of the cardinality of nonstandard finite sets and factor it out by the convex hull of the nonstandard reals. In the later papers \cite{Hru, Hru2}, Hrushovski states some properties of this pseudofinite dimension and used it to get asymptotic results in additive  combinatorics.

We present a similar connection between forking and the logarithmic pseudofinite dimension : any instance of dividing in a pseudofinite structure is witnessed by a decrease of the dimension. This connection is used to get some known results in asymptotic classes of finite structures, as defined in \cite{MS} and \cite{Elw}.

The paper is organized as follows: in Section \ref{quasidimensions} we present the definition of quasidimensions and the construction of the logarithmic pseudofinite dimension. In Section \ref{pseudoforking} we present the main result of this paper: any instance of dividing can be witnessed by a decrease of the pseudofinite dimension. We start the proof recalling some results in combinatorics and measure theory (Section \ref{measure}) and using them to prove the result in Section \ref{proofpseudoforking}.

Section \ref{asympclass} contains a calculation of the possible pseudofinite dimensions for the ultraproducts of 1-dimensional asymptotic classes, which can be easily generalized to the context of $N$-dimensional asymptotic classes. As a corollary, we obtain new proofs of the following known results from \cite{MS},\cite{Elw}: Every infinite ultraproduct of the members of a $1$-dimensional asymptotic class (resp. $N$-dimensional asymptotic class) is supersimple of $U$-rank 1. (resp. $U$-rank less than or equal to $N$).\\

\noindent \textbf{Acknoledgements:} I would like to thank my advisors Thomas Scanlon and Alf Onshuus for all their help and their valuable comments and discussions through the development of this paper.

\section{The logarithmic pseudofinite dimension} \label{quasidimensions}
In this section we present the definition of quasidimension as presented in \cite{HW} and give the construction of the \emph{logarithmic pseudofinite dimension}, also proving that it define a quasidimension on the ultraproducts of finite structures.
 
\bd Let $M$ be any structure. A \emph{quasi-dimension} on $M$ is a map $\delta$ from the class of definable sets into an ordered abelian group $G$, together with a formal element $-\infty$, satisfying:
\benum
\item $\delta(\emptyset)=-\infty$, and $\delta(X)>-\infty$ implies $\delta(X)\geq 0$.
\item $\delta(X\cup Y)=\max\,\{\delta(X),\delta(Y)\}$
\item For every $g\in G\cup\{-\infty\}$ the following holds: If $X$ is a definable subset of $M^k$, $\pi$ is the projection to some of the coordinates and  $\delta(\pi^{-1}(\ov{x}))\leq g$ for all $\ov{x}\in \pi(X)$, then $\delta(X)\leq \delta(\pi(X))+g$.
\eenum
\ed 

We will focus in the logarithmic pseudofinite dimension, which is a quasidimension defined on ultraproducts of finite structures. Consider the following construction:

Assume $M$ is an infinite ultraproduct of finite structures $\langle M_i:i\in I\rangle$, with $\displaystyle{|M_i|\to \infty}$. For a definable set $X$, there is a map

\begin{align*}
\log_i: &\operatorname{Def}(M_i)\longrightarrow \R \cup \{-\infty\} \\
&\ \ X(M_i)\, \longmapsto \log(|X(M_i)|)
\end{align*}
where $\log$ is the usual natural logarithm and $|X(M_i)|$ represents the size of the definable set $X(M_i)$.

It is possible to take the ultraproduct of such functions and obtain a map

\begin{align*}
\log = \prod_{\mathcal{U}} \log_i: &\operatorname{Def}(M)\longrightarrow \R^*\cup \{-\infty\} \\
&\ \ \ \ X \ \ \ \  \longmapsto \log(|M|):=[\log(|X(M_i)|)]_{i\in I}
\end{align*}
where $\R^*$ is a non-standard real closed field.
Let $\mathcal{C}$ be the convex hull of $\Z$ in $\R^*$ (a convex subgroup of $\R^*$) and $\pi:\R^* \cup\{-\infty\}\rightarrow \left(\R^*/\mathcal{C}\right) \cup \{-\infty\}$ the natural quotient map (with $\pi(-\infty)=-\infty$). 

For a definable subset $X$ of $M$, define \[\delta(X)=\pi(log(|X|))\]

This is a way to measure ``bigness'' of the definable sets in $M$. For instance, note that $\delta(X)=0$ if and only if $\log(|X|)\in \mathcal{C}$, which (by compactness) implies that $|X_i|$ is uniformly bounded by a fixed $M$ on a $\mathcal{U}$-large set.
\bp $\delta(X)=\pi(log(|X|))$ is a quasi-dimension on $M=\prod_{\mathcal{U}} M_i$
\ep
\bdem
\benum
\item[(1)] Clearly $\delta(\emptyset)=-\infty$, and for any definable $X$ we have that, if $X$ is non-empty in the ultraproduct, then \[\{i\in I:|X_i|\geq 1\}=\{i\in I:\,log_i(|X_i|)\geq 0\}\in \mathcal{U}\] which implies $\delta(X)\geq 0$.
\item[(2)] Let $X,Y$ be definable subsets of $M$, and assume without loss of generality that $\delta(X)\geq \delta(Y)$. By the construction of $\delta$, this implies in particular that $|X(M_i)|\geq |Y(M_i)|$ for $\mathcal{U}$-almost all $i\in I$. In those indices, we have
\begin{align*}
|X(M_i)| &\leq |X(M_i)\cup Y(M_i)|\leq 2|X(M_i)|\\ 
\log(|X(M_i)|) &\leq \log(|X(M_i)\cup Y(M_i)|)\leq \log(2)+\log(|X(M_i)|)\\
0 &\leq \log(|X(M_i)\cup Y(M_i)|)-\log(|X(M_i)|)\leq \log(2)
\end{align*}
So, $\delta(X\cup Y)=\pi(\log(|X\cup Y|))=\pi(\log(|X|))=\delta(X)$.
\item[(3)] Let $X\subseteq^{def} M^k$ and $p$ be a projection to some of the coordinates and $\alpha\in \R^*/\mathcal{C}$. Also assume $\delta(p^{-1}(\ov{x}))\leq \alpha$ for all $\ov{x}\in p(X)$. 

Since $X=\prod_\mathcal{U} X(M_i)$ we have $\mathcal{U}$-almost everywhere that:
\begin{align*}
|X(M_i)|&\leq |p(X(M_i))|\cdot \max\,\{|p^{-1}(\ov{x})|: \ov{x}\in p(X(M_i))\} \\
\log(|X(M_i)|)&\leq \log(|p(X(M_i))|) + \max\, \{\log(|p^{-1}(\ov{x})|):\ov{x}\in p(X(M_i))\}
\end{align*}
and since the natural projection $\pi:\R^*\longrightarrow \R^*/\mathcal{C}$ is an order-preserving group homomorphism, we conclude that

\[\delta(X)\leq \delta(p(X))+\sup\,\{\delta(p^{-1}(\ov{x})):\ov{x}\in p(X)\}\leq \delta(p(X))+\alpha.\]
\eenum
\edem

\section{Pseudofinite dimension and forking} \label{pseudoforking}

The purpose of this section is to show the relationship between the pseudofinite dimension defined in Section \ref{quasidimensions} and the forking relation inside the structure $M$. This relationship can be viewed as a generalization of the notion of rank in stable or simple theories, in the sense that dividing can be witnessed by a drop in the dimension. 

To prove this, we will use some more or less known results in enumerative combinatorics and measure theory. We include the proofs here for completeness.

\subsection{Some lemmas from combinatorics and measure theory} \label{measure}

We start with the following lemma:

\bl If $m\geq 2i+1$ are integers, then $\displaystyle{\binom{m}{i} - \binom{m}{i+1}\leq 0}$
\el
\bdem This is a straightforward calculation:
\begin{align*}
\binom{m}{i}-\binom{m}{i+1}&=\dfrac{m!}{i!(m-i)!}-\dfrac{m!}{(i+1)!(m-i-1)!}\\
&=\dfrac{m!}{i!(m-i-1)!}\left(\dfrac{1}{m-i}-\dfrac{1}{i+1}\right)\\
&=\dfrac{m!}{i!(m-i-1)!}\left(\dfrac{2i+1-m}{(m-i)(i+1)}\right)\\
&\leq \dfrac{m!}{(i+1)!(m-i)!}(2i+1-2i-1)=0.
\end{align*}
\edem

Now we present the measure theoretic lemma. Assume we have a measure space $(X,\mathcal{B},\mu)$. Given measurable sets $A_1,\ldots,A_n$, we can define $S_i$ to be the sum of the measures of all $k$-intersections of $A_1,\ldots,A_n$, namely,

\[S_k:=\sum_{1\leq i_1<i_2<\cdots<i_k\leq n} \mu(A_{i_1}\cap \cdots \cap A_{i_k})\]

We know from the inclusion-exclusion principle that the measure of $\displaystyle{\bigcup_{i=1}^n A_i}$ is an alternating sum of $S_i's$. What we will prove now is that starting from positive term of this sum, the result is positive. 

\bp \label{truncatedinclusionexclusion} \emph{[Truncated inclusion-exclusion principle]}
Let $X$ be a measure space and $A_1,\ldots,A_n$ be measurable sets, and let $S_1,\ldots,S_n$ as defined above. Then for every $k\leq n/2$, 
\[\sum_{i=2k+1}^n (-1)^{i-1}S_i\geq 0\]
\ep
\bdem By the inclusion-exclusion principle we know that 
\[\mu\left(\bigcup_{i=1}^n A_i\right)=\sum_{m=1}^{n} (-1)^{m-1}S_m\]

For every non-empty $W\subseteq \{1,\ldots,n\}$, define $\displaystyle{E_W=\bigcap_{i\in W} A_i\cap \bigcap_{i\not\in W}A_i^c.}$
These are the non-empty atoms of the algebra of sets generated by $A_1,\ldots,A_n$ that are contained in $\displaystyle{\bigcup_{i=1}^n A_i}$. They are disjoint and we have the following easy identities:
\[\mu\left(\bigcup_{i=1}^n A_i\right)=\sum_{W\subseteq \{1,\ldots,n\}} \mu(E_W), \hspace{1cm}\mu(A_{i_1}\cap \cdots \cap A_{i_m})=\sum_{i_1,\ldots,i_m\in W}\mu(E_W)\]

So, the inclusion-exclusion principle states that
\begin{align*}
\sum_{W\subseteq \{1,\ldots,n\}} \mu(E_W) &:=\sum_{m=1}^n (-1)^{m-1} S_m\\
&=\sum_{m=1}^n (-1)^{m-1} \left(\sum_{1\leq i_1<i_2<\cdots<i_m\leq n} \mu(A_{i_1}\cap \cdots \cap A_{i_m})\right)\\
&=\sum_{m=1}^n (-1)^{m-1} \left(\sum_{1\leq i_1<i_2<\cdots<i_m\leq n} \left(\sum_{i_1,\ldots,i_m\in W} \mu(E_W)\right)\right)\\
&=\sum_{m=1}^n (-1)^{m-1} \left(\sum_{W\subseteq \{1,\ldots,n\}}\alpha^m_W\mu(E_W)\right)
\end{align*}
where the coefficient $\alpha^m_W$ is the number of times that the summand $\mu(E_W)$ appears in $S_m$.
Thus, we know that for every $W\subseteq \{1,\ldots,n\}$, $\displaystyle{\sum_{m=1}^n(-1)^{m-1}\cdot \alpha^m_W=1}$

So we have three cases:
\bitem
\item If $|W|<2k+1$, all the possible summands $\mu(E_W)$ are already in the sum $\displaystyle{\sum_{i=1}^{2k} (-1)^{i-1}S_i.}$
\item If $|W|=2k+1$, the coefficient of $\mu(E_W)$ in $\displaystyle{\sum_{i=1}^{2k} (-1)^{i-1}S_i}$ is $0$, because there is exactly one more term $\mu(E_W)$ appearing in $S_{2k+1}$ (the measure of the intersection $\displaystyle{\bigcap_{i\in W} A_i}$)
\item If $|W|>2k+1$: Note that in every $S_i$ appear $\binom{|W|}{i}$ summands of the form $\mu(E_W)$ for a fixed $W\subseteq\{1,\ldots,n\}$. So, the coefficient of $\mu(E_W)$ in $\displaystyle{\sum_{i=1}^{2k} (-1)^{i-1}S_i}$ is:
\begin{align*}
\beta_W=\sum_{m=1}^{2k} (-1)^{m-1}\binom{|W|}{i}=\sum_{j=1}^k \left[ \binom{|W|}{2j-1}-\binom{|W|}{2j}\right]\leq 0
\end{align*} using the previous lemma (note that $j\leq k$ implies $2j-1\leq 2k-1<|W|$).
\eitem
Therefore, we obtain
\begin{align*}
\sum_{W\subseteq\{1,\ldots,n\}}\mu(E_W)&=\sum_{m=1}^n(-1)^{m-1}S_m\\
&=\sum_{m=1}^{2k} (-1)^{m-1}S_m + \sum_{m=2k+1}^n (-1)^{m-1}S_m\\
&=\sum_{|W|<2k+1}\mu(E_W)+0+\left(\sum_{|W|=2k+1}\beta_W\cdot \mu(E_W)\right) + \sum_{m=2k+1}^n (-1)^{m-1}S_m\\
\end{align*}

So, \[\sum_{|W|>2k+1}\mu(E_W)=\sum_{|W|=2k+1}\beta_W\cdot \mu(E_W) + \sum_{m=2k+1}^n (-1)^{m-1}S_m,\]
and we conclude that
\[\sum_{m=2k+1}^n (-1)^{m-1}S_m=\sum_{|W|>2k+1}\mu(E_W)-\sum_{|W|=2k+1}\beta_W\mu(E_W)\geq 0\]because all the coefficients $\beta_W$ are less than or equal to $0$.
\edem

The following measure-theoretic proposition will play a key role in the proof that dividing implies a decrease of pseudofinite dimension.

\bp \label{kintersections} Let $X$ be a measure space with $\mu(X)=1$ and fix $0<\epsilon\leq \dfrac{1}{2}$. Let $\langle A_i:i<\omega\rangle$ be a sequence of measurable subsets of $X$ such that $\mu(A_i)\geq \epsilon$ for every $i$

Then, for every $k<\omega$ there are $i_1<i_2<\ldots<i_k$ such that \[\mu\left(\bigcap_{j=1}^k A_{i_j}\right)\geq \epsilon^{3^{k-1}}\] 
\ep
\bdem The proof will be by induction on $k$. 
\bitem
\item $k=1$: By hypothesis we have $\mu(A_i)\geq \epsilon= \epsilon^{3^{1-1}}$.  \checkmark
\item $k=2$: Assume not, then $\mu(A_i\cap A_j)<\epsilon^{3^{2-1}}=\epsilon^3$ for every $i\neq j$. By the truncated inclusion-exclusion principle we know that for every $N\in\Nat$,
\begin{align*}
1 \geq \mu\left(\bigcup_{i=1}^N A_i\right) &\geq \sum_{i=1}^N \mu(A_i) - \sum_{1\leq i<j\leq N} \mu(A_i\cap A_j) &(\text{by Proposition \ref{truncatedinclusionexclusion}})\\
&\geq N\epsilon - \dfrac{N(N-1)}{2}\epsilon^3 &(\dag)
\end{align*}
Define the quadratic function given by \[f(x)=x\cdot \epsilon - \dfrac{x(x-1)}{2}\epsilon^3=-\dfrac{x^2}{2}\epsilon^3 + x\left(\epsilon+\dfrac{\epsilon^3}{2}\right).\]

This function achieve its maximum value at $x_0=\dfrac{1}{\epsilon^2}+\dfrac{1}{2}$>0, and by taking any integer $N\in [x_0-1,x_0]$ we have that 
\begin{align*}
f(N)\geq f(x_0-1)&=\left(\dfrac{1}{\epsilon^2}-\dfrac{1}{2}\right)\epsilon - \dfrac{\left(\dfrac{1}{\epsilon^2}-\dfrac{1}{2}\right)\left(\dfrac{1}{\epsilon^2}-\dfrac{3}{2}\right)}{2} \cdot \epsilon^3\\
&=\dfrac{1}{\epsilon}-\dfrac{\epsilon}{2}- \dfrac{\dfrac{1}{\epsilon^4}-\dfrac{2}{\epsilon^2}+\dfrac{3}{4}}{2}\cdot \epsilon^3\\
&=\dfrac{1}{\epsilon}-\dfrac{\epsilon}{2}-\dfrac{1}{2\epsilon}+\epsilon-\dfrac{3}{8}\epsilon^3\\
&=\dfrac{1}{2\epsilon}+\dfrac{\epsilon}{2}-\dfrac{3}{8}\epsilon^3\\
&\geq 1+\epsilon \left(\dfrac{1}{2}-\dfrac{3}{8}\epsilon^2\right) \hspace{1cm}\text{[because $\epsilon\leq \frac{1}{2}$]}\\
&>1.
\end{align*} contradicting the inequality (\dag). \checkmark
\eitem
Now, assume the induction hypothesis, which is that there is a tuple $(i_1,\ldots,i_k)$ satisfying \[i_1<\ldots < i_k \text{\,\,\, and \,\,\,} \mu\left(\bigcap_{j=1}^k A_{i_j}\right)\geq \epsilon^{3^{k-1}}\text{\hspace{1cm}} (*)\]

\emph{Claim: There are infinitely many such tuples}.

\emph{Proof of the Claim:} Assume not, and take $\ell$ to be the maximum of all indices appearing in the tuples $(i_1,\ldots,i_k)$ which satisfies (*). The sequence $\langle A_{j}:j\geq \ell+1\rangle$ would contradict the induction hypothesis. \checkmark \\

Now, let $\langle \alpha_j:j<\omega\rangle$ be an enumeration of all tuples satisfying (*) and define $\displaystyle{B_j=\bigcap_{i\in \alpha_j}A_i}$. By construction, $\mu(B_j)\geq \epsilon^{3^{k-1}}$.

By the $k=2$ case, there are indices $j_1\neq j_2$ such that \[\mu(B_{j_1}\cap B_{j_2})\geq (\epsilon^{3^{k-1}})^3=\epsilon^{3^{k-1}\cdot 3}=\epsilon^{3^k}\] where $j_1,j_2$ are indices corresponding to two different tuples $\alpha_{j_1}\neq \alpha_{j_2}$. In particular, there are (at least) $k+1$ indices $i_1<i_2<\cdots<i_k<i_{k+1}$ such that
\[\mu\left(\bigcap_{j=1}^{k+1}A_{i_j}\right)\geq \mu(B_{j_1}\cap B_{j_2})\geq \epsilon^{3^k}=\epsilon^{3^{(k+1)-1}} \hspace{1cm}\]  
\edem

\subsection{Dividing and drop of the pseudofinite dimension}\label{proofpseudoforking}

With the results of the previous subsection, we are now able to give a proof of the main result of this note. The setting, as in the definition of logarithmic pseudofinite dimension, is the following: $\langle M_i:i\in I\rangle$ is a family of finite structures, $M$ is an infinite ultraproduct of the family and $\delta$ denotes the logarithmic pseudofinite dimension defined on definable subsets of $M$. 

\bt \label{forkingdelta} Let $X=\psi(x,\ov{a})$ be a definable subset of $M$ and $\phi(x,\ov{b})$ a formula implying $\psi(x,\ov{a})$. If $\phi(x,\ov{b})$ divides over $\ov{a}$, then there exists $\ov{b}'\models tp(\ov{b}/\ov{a})$ such that $\delta(\phi(x,\ov{b}'))<\delta(X)$.
\et

\bdem Assume that $\phi(x,\ov{b})$ divides over $\ov{a}$. Towards a contradiction, assume that for every $b'\models tp(\ov{b}/\ov{a})$ we have $\delta(\phi(x,b'))=\delta(X)$. Then for each $\ov{b}'$ there is $n_{b'}\in\Nat$ such that \[\log(|X|)-\log(|\phi(x,b')|)<n_{b'}.\] 

Thus, 
\begin{align*}
\log\left(\dfrac{|X|}{|\phi(x,b')|}\right)&<n_{b'}\\
\dfrac{|X|}{|\phi(x,b')|}&<e^{n_{b'}}\\
\dfrac{|\phi(x,b')|}{|X|}&\geq e^{-n_{b'}}\\
e^{n_{b'}}|\phi(x,b')|&\geq |X|
\end{align*}
In particular, there is $M_{b'}\in\Nat$ such that $M_{b'}|\phi(x,b)|\geq |X|$.\\

\emph{\textbf{Claim:} There is an uniform bound $M$ such that $M|\phi(x,b')|\geq |X|$ for every $b'\models tp(b/\ov{a})$.}\\

\noindent \emph{Proof of the Claim:} If not, for every $n<\omega$ there is $\ov{b}_n\equiv_{\ov{a}} \ov{b}$ such that $log( |X|)-log(|\phi(x,\ov{b}_n)|)>n$. Consider the multi-sorted structures given by \[\mathcal{M}_i=\langle M_i,\R,log_{\varphi}\rangle_{\varphi\in\Le}\] where $log_\varphi$ is an function between different sorts interpreted as 
\begin{align*}
log_\varphi: &M_i\longrightarrow \R\\
&\ov{b} \longmapsto \log(|\varphi(x,\ov{b})|)
\end{align*} 

Now, take the ultraproduct $\mathcal{M}=\prod_{\mathcal{U}} \mathcal{M}_i$. In this structure, consider the type \[p(y)=tp(\ov{b}/\ov{a})\cup \{log_{\psi}(\ov{a})-log_{\phi}(\ov{b})>n:n\in \Nat\}\] This type is finitely satisfiable in $\mathcal{M}$, and by $\aleph_1$-saturation of the ultraproduct, there is $\ov{b}'\models p(y)$ which means $\ov{b}'\models tp(\ov{b}/\ov{a})$ and $\delta(\phi(x,\ov{b}'))<\delta(X)$, a contradiction. \hspace{0,5cm}\checkmark \\

Therefore, there is $M\in\Nat$ such that $M\cdot |\phi(x,\ov{b}')|\geq |X|$ for every $\ov{b}'\equiv_{\ov{a}} \ov{b}$. Since $\phi(x,\ov{b})$ divides over $\ov{a}$, there is an indiscernible sequence $\langle \ov{b}_j:j<\omega\rangle$ (which can be assumed to be in $M$ by $\aleph_1$-saturation) such that:
\bitem
\item $\ov{b}_i\models tp(\ov{b}/A)$.
\item $\{\phi(x,\ov{b}_j):j<\omega\}$ is $k$-inconsistent for some $k<\omega$.
\eitem

Assume $\ov{b}_j=[\ov{b}_j^{i}]_{\sim_{\mathcal{U}}}$. By the claim, $M\cdot |\phi(x,\ov{b}_j)|\geq |X|$ which  implies that $\displaystyle{\frac{|\phi(x,\ov{b}_j^{i})|}{|X|}\geq \frac{1}{M}}$ for $\mathcal{U}$-almost all $i$. 

Consider the normalized counting measure localized on $X(M_i)$ in each finite structure $M_i$, and the Loeb measure induced on $M$ by these measures. Therefore, $\mu(X)=1$ and $\langle \phi(M,\ov{b}_j):j<\omega\rangle$ is a sequence of measurable sets with $\mu(\phi(M,\ov{b}_j))\geq \frac{1}{M}$ for every $j<\omega$. By the Proposition \ref{kintersections} there are $j_1<\ldots <j_k<\omega$ such that \[\mu\left(\bigcap_{l=1}^k \phi(M,\ov{b}_{j_l})\right)\geq \dfrac{1}{M^{3^{k-1}}}>0,\]in particular, $\displaystyle{\bigcap_{l=1}^k \phi(M,b_{j_l})}$ is non-empty, contradicting $k$-inconsistency.
\edem

The theorem above allows us to conclude that the number of possible different values for pseudofinite dimensions of definable sets is a bound for the length of dividing chains, providing also a bound for the U-rank in types. We will explore this idea in the next section.\\

We might think about two possible generalizations of Theorem \ref{forkingdelta}: either changing dividing by forking or showing that the original formula has lower pseudofinite dimension. The following two examples answer negatively to these attempts:

\bej Consider the class of finite structures $M_n=([0,2^n],E_n)$ where $E_n$ is an equivalence relation with classes \[[0,2^{n-1}-1],[2^{n-1},2^{n-1}+2^{n-2}-1],[2^{n-1}+2^{n-2},2^{n-1}+2^{n-2}+2^{n-3}-1],\ldots,[2^{n-1}+2^{n-2}+\cdots +2^{2},2^n]\] Let $M=\prod_\mathcal{U} M_n$ and $b=[\ov{0}]_{\mathcal{U}}$. In the ultraproduct  $M$ the relation symbol $E$ is interpreted as an equivalence relation with infinitely many infinite classes, and as it was shown in Example \ref{equivalencerelationforking} we have that the formula $xEb$ divides over the empty set. Theorem \ref{forkingdelta} shows that there is a conjugate of $b$ over $\emptyset$ witnessing a drop of pseudofinite dimension, but this drop is not witnessed by the formula $xEb$ because

\[\log |M_n|-\log |xE_n0|=\log (2^n)- \log (2^{n-1})=\log 2<1\] which implies that $\delta(M)=\delta(xEb)$.
\eej

\bej This example is an adaptation of the classical example of the circular order that shows that the formula $x=x$ may fork over $\emptyset$. Consider the structures $M_n=(\Z/(3n)\Z,R)$ where $R$ is a ternary relation interpreted in $M_n$ as follows: $M_n\models R(b;\ov{a},\ov{c})$ if and only if there are integers $a',b',c'$ congruent with $a,b,c$ $(\mod\,3n)$ respectively, such that $a'<b'<c'$ and $|c'-a'|<n+1$. 

Take $M=\prod_\mathcal{U}M_n$, and the elements $a:=[a_n=0]_{\mathcal{U}},b:=[b_n=n]_{\mathcal{U}}\in M$. \\

\noindent \emph{Claim:} The formula $R(x;a,b)$ divides over $\emptyset$. \\

\noindent \emph{Proof of the Claim:} On each $M_n$ we can define the tuples Consider the sequence given by $\left\langle (a^n_i,b^n_i)=(n+k\cdot \llbracket\log n\rrbracket,n+(k+1)\cdot \llbracket\log n\rrbracket):k\leq \frac{n}{\log n}\right\rangle$, and consider in the ultraproduct the sequence given by $\left\langle (a_i,b_i)=([a^n_i]_{n\in\mathcal{U}},[b^n_i]_{n\in\mathcal{U}}):i<\omega\right\rangle$. This is a sequence in $tp(a,b/\emptyset)$ which is indiscernible over the empty set, and by construction we have that the set of formulas $\{R(x;a_i,b_i):i<\omega\}$ is 2-inconsistent. \hspace{0.5cm}\checkmark\\

Consider the elements in the ultraproducts given by 
\[a_1:=[a_1^n=0]_{\mathcal{U}}=b_3,a_2:=[a_2^n=n]_{\mathcal{U}}=b_1 \text{\ \ and\ \ }a_3:=[a_3^n=2n]_{\mathcal{U}}=b_2.\]

Note that the formula $x=x$ forks over $\emptyset$, because it implies the disjunction \[\bigvee_{i=1}^3 R(a_i,x,b_i)\vee \bigvee_{i=1}^3 x=a_i\]which is a disjunction of formulas that divide over $\emptyset$.\\

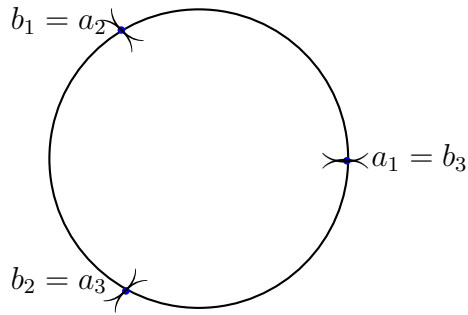
\begin{center}
\begin{pspicture}(0,0)(4,4)
\pscircle(2,2){2}
\put(4.3,1.9){$a_1=b_3$}
\put(-0.5,3.728){$b_1=a_2$}
\put(-0.5,0,272){$b_2=a_3$}

\psdots[dotsize=0.1cm,linecolor=blue](3.97,1.97)
\psdots[dotsize=0.1cm,linecolor=blue](0.97,3.712050)
\psdots[dotsize=0.1cm,linecolor=blue](1.03,0.23794)
\rput{90}(3.95,2){\Large{(}}
\rput{90}(3.95,1.95){\Large{)}}

\rput{45}(1,3.732050){\Large{(}}
\rput{45}(0.95,3.682050){\Large{)}}

\rput{135}(1,0.26794){\Large{(}}
\rput{135}(1.05,0.21794){\Large{)}}

\end{pspicture}
\end{center}

However, the realization of the formula of $x=x$ is $M$ and it does not witness any drop of pseudofinite dimension ($\delta(M)$ is the maximal value of the pseudofinite dimension among subsets of $M$).
\eej

\section{Pseudofinite dimension and asymptotic classes} \label{asympclass}

In general, the logarithmic pseudofinite dimension can take infinitely many different values on the definable sets of $M=\prod_{\mathcal{U}} M_i$. For instance, consider the class $\mathcal{C}_{ord}$ of finite linear orders. If $M_n=([1,n],<)$ and $\alpha,\beta$ are elements in the interval $[0,1]$ with $\alpha<\beta$, we may define $X(M_n)=[1,\llbracket n^\alpha\rrbracket]$ and $Y(M^n)=[1,\llbracket n^\beta\rrbracket]$. We will show that $\delta(X)<\delta(Y)$.

Otherwise, $\delta(X)$ and $\delta(Y)$ are equal, which means there is a natural number $N$ such that $\log |X|+N\geq \log |Y|.$ Therefore, 
\begin{align*}
\alpha \log n + N \geq \beta\log n\\
N\geq (\beta-\alpha)\log n
\end{align*}which is not true when $n$ tends to infinity.

The main feature of this example is that an infinite linear order can be defined on the ultraproducts, implying they are unstable non-simple and thus they have arbitrarily long dividing chains.

On the other hand, there are classes of finite structures with a better behavior of their ultraproducts. That is the case of the 1-dimensional asymptotic classes (and more generally of the N-dimensional asymptotic classes) whose definition appear in \cite{MS} and \cite{Elw}. These classes are known to have supersimple ultraproducts, which implies a finite bound on the length of dividing chains in their ultraproducts.\\

The purpose of this section is to show that the supersimplicity of these classes can be detected by the logarithmic pseudofinite dimension. For instance, we will show that for 1-dimensional classes (which ultraproducts are supersimple of U-rank $1$) the only possible values for the pseudofinite dimension are $-\infty,0$ and $\alpha=\delta(M)$.\\

First, recall the definition of these classes:

\bd \label{onedimclass} Let $\Le$ be a first order language, and $\C$ be a collection of finite $\Le$-structures. Then $\C$ is a \emph{1-dimensional asymptotic class} if the following hold for every $m\in\Nat$ and every formula $\varphi(x,\ov{y})$, where $\ov{y}=(y_1,\ldots,y_m)$.
\benum
\item There is a positive constant $C$ and a finite set $E\subseteq \R^{>0}$ such that for every $M\in\C$ and $\ov{a}\in M^m$, either $|\varphi(M,\ov{a})|\leq C$ or for some $\mu\in E$,
\[||\varphi(M,\ov{a})|-\mu |M||\leq C|M|^{1/2}.\]
\item For every $\mu\in E$, there is an $\Le$-formula $\varphi_\mu(\ov{y})$ such that, for all $M\in\C$, $\varphi_\mu(M^m)$ is precisely the set of $\ov{a}\in M^m$ with 
\[||\varphi(M,\ov{a})|-\mu |M||\leq C|M|^{1/2}.\]
\eenum
\ed

\bp \label{delta01} Let $\mathcal{C}$ be a class of finite structures. If $\mathcal{C}$ satisfies the condition (1) in  definition \ref{onedimclass}, then for every infinite ultraproduct $M$ of elements in $\mathcal{C}$ there are only two possible values for $\delta(X)$ while $X$ varies among the non-empty definable subsets of $M^1$.
\ep
\bdem Let $\varphi(x,\ov{a})$ be a definable set in the ultraproduct and take $\mu_1,\ldots,\mu_k>0$ the possible measures in $E$ satisfying (1). Assume $M=\prod_{\mathcal{U}}M_i$ with $i\in I$ and $\ov{a}=[\ov{a_i}]_{i\in I}$.

For every $M_i \in \mathcal{C}$ one of the following hold:
\bitem
\item $|\varphi(M_i,\ov{a_i})|\leq C$
\item $\left| |\varphi(M_i,\ov{a_i})| - \mu_j |M_i|\right|\leq C|M_i|^{1/2}$ for some $j=1,2,\ldots,k$.
\eitem 

Consider the sets
\begin{align*}
A_0&=\{i\in I: |\varphi(M_i,\ov{a_i})|\leq C\}\\
A_1&=\{i\in I: \left| |\varphi(M_i,\ov{a_i})| - \mu_1 |M_i|\right|\leq C|M_i|^{1/2}\}\\
&\vdots \\
A_k&=\{i\in I: \left| |\varphi(M_i,\ov{a_i})| - \mu_k |M_i|\right|\leq C|M_i|^{1/2}\}
\end{align*}
Since $A_0\cup A_1 \cup \cdots \cup A_k=I$, one of these sets belongs to $\mathcal{U}$ because $\mathcal{U}$ is an ultrafilter on $I$. We have to consider two cases:

\bitem
\item If $A_0\in\mathcal{U}$, then $|\varphi(M_i,\ov{a_i}|\leq C$ (a.e. in $\mathcal{U}$) for some fixed $C>0$, which implies $|\varphi(M,\ov{a})|\leq C$ and therefore \[\delta(\varphi(M,\ov{a}))=\pi(\log(|\varphi(x,\ov{a})|))\leq \pi(\log(C))=0\]

\item If $A_j\in\mathcal{U}$ for some $j=1,2,\ldots,k$ then we obtain
\[\mu_j|M_i|-C|M_i|^{1/2} \leq |\varphi(M_i,\ov{a_i})|\leq \mu_j|M_i|+C|M_i|^{1/2}\]

Put $\mu_*=\min\{\mu_1,\ldots,\mu_k\}$ and $\mu^*=\max\{\mu_1,\ldots,\mu_k\}$. Then for every definable set $X=\prod_\mathcal{U} X_i$, either $\delta(X)=0$ (because the corresponding set $A_0$ belong to $\mathcal{U}$) or 
\[\mu_* |M_i|-C|M_i|^{1/2} \leq |X_i|\leq \mu^*|M_i|+C|M_i|^{1/2}\]
So, 
\begin{align*}
|M_i|^{1/2}\left(\mu_* |M_i|^{1/2} - C\right) &\leq |X_i| \leq |M_i|^{1/2}\left(\mu^* |M_i|^{1/2} + C\right) &\\
|M_i|^{1/2}\left(\mu_* |M_i|^{1/2} - \dfrac{\mu_*}{2}|M_i|^{1/2} \right)& \leq |X_i| \leq |M_i|^{1/2}\left(\mu^* |M_i|^{1/2} + \dfrac{\mu^*}{2}|M_i|^{1/2}\right) \\
& \text{(asymptotically, because $|M_i|\to\infty$)}\\
|M_i|^{1/2}\left(\dfrac{\mu_*}{2}|M_i|^{1/2} \right) & \leq |X_i| \leq |M_i|^{1/2}\left( \dfrac{3\mu^*}{2}|M_i|^{1/2}\right) \\
\dfrac{1}{2}\log(|M_i|)+\log\left(\dfrac{\mu_*}{2}\right)+\dfrac{1}{2}\log(|M_i|) &\leq \log(|X_i|)\leq \dfrac{1}{2}\log(|M_i|)+\log\left(\dfrac{3\mu^*}{2}\right)+\dfrac{1}{2}\log(|M_i|)\\
\log(|M_i|)+\log\left(\dfrac{\mu_*}{2}\right)&\leq \log(|X_i|)\leq \log(|M_i|)+\log\left(\dfrac{3\mu^*}{2}\right)\\
\pi(\log(|M_i|))=\pi\left(\log(|M_i|)+\log\left(\dfrac{\mu_*}{2}\right)\right) &\leq \pi(\log(|X_i|))\\
&\leq \pi\left(\log(|M_i|)+\log\left(\dfrac{3\mu^*}{2}\right)\right)=\pi(\log(|M_i|))\\
\delta(M)&\leq \delta(X)\leq \delta(M)
\end{align*}as we desired.
\eitem
\edem
Now, we present a new proof of the following known result, that appears in \cite{MS} (Lemma 4.1).
\bc [cf. \cite{MS}, Theorem 4.1] Let $\mathcal{C}$ be a class of finite structures satisfying the condition (1) of Definition \ref{onedimclass}, and let $M$ be an infinite ultraproduct of members of $\mathcal{C}$. Then $Th(M)$ is supersimple of $SU$-rank 1.
\ec
\bdem 
Assume $SU(M)\geq 2$. Then there is an increasing chain of types $p_0\subset p_1\subset p_2$ such that $p_{i+1}$ is a dividing extension of $p_i$ for $i=0,1$. In particular, there are formulas $\phi_1(x,a_1)\in p_1$ which divides over $A_0$ and $\phi_2(x,a_2)\in p_2$ which divides over $A_1$. By Proposition \ref{forkingdelta} there are tuples $a_1'\equiv_{A_0} a_1$ and $a_2'\equiv_{A_1}a_2$ such that \[\delta(\phi(x,a_2'))<\delta(\phi(x,a_1'))<\delta(M)\] This contradicts Proposition \ref{delta01} which states there are only two possible values for $\delta$ on non-empty definable subsets of $M$.
\edem




\begin{thebibliography}{20}
\bibitem{ChHru} G. Cherlin. E. Hrushovski. \emph{Finite strucutres with few types}. Princeton University Press. Princeton and Oxford. 2003
\bibitem{Elw} R. Elwes. \emph{Asymptotic classes of finite structures.} Journal of Symbolic Logic. Volume 72, Issue 2 (2007), 418-438.
\bibitem{GolTwo} I. Goldbring. H. Towsner. \emph{An approximate logic for measures}. Israel Jornal of Mathematics (to appear). Preprint. \url{arXiv:1106.2854v1}. June 2011.
\bibitem{Hodges} W. Hodges. \emph{Model Theory}. Encyclopedia of mathematics and its applications. Cambridge University Press. 1994
\bibitem{Hru} E. Hrushovski. \emph{Stable group theory and approximate subgroups}. Journal of the American Mathematical Society. Volume 25. Number 1. January 2012. Pages 189-243.
\bibitem{Hru2} E. Hrushovski. \emph{On Pseudo-Finite Dimensions.} Notre Dame Journal of Formal Logic. Volume 54 (2013), no. 3-4, 463--495.
\bibitem{HW} E. Hrushovski, F. Wagner. \emph{Counting and dimensions}. Model Theory with applications to Algebra and Analysis. 2008
\bibitem{MS} D. Macpherson, C. Steinhorn. \emph{One-dimensional asymptotic classes of finite structures.} Transactions of the American Mathematical Society. Volume 360, pages 411-448. 2007.
\bibitem{TaoVu} T. Tao. V.H. Vu \emph{Additive Combinatorics}. Cambridge studies in advances mathematics [105]. Cambridge University Press. 2006
\end{thebibliography}
\end{document}